\documentclass[11pt]{amsart}
\usepackage[margin=2.35cm]{geometry}
\usepackage{amssymb,mathtools}
\usepackage{booktabs}
\usepackage{enumitem}
\usepackage[T1]{fontenc}
\usepackage{lmodern}
\usepackage{microtype}
\usepackage{xcolor}
\usepackage[colorlinks=true,linkcolor=blue!50!black,citecolor=blue!50!black,urlcolor=blue!50!black,bookmarks=true,bookmarksnumbered=true,unicode=true]{hyperref}

\hypersetup{pdftitle={More lower bounds for Hales--Jewett and van der Waerden numbers},pdfauthor={Younes Mouhib},pdfsubject={Explicit colourings, algebraic certificates, local search, and exact SAT methods},pdfkeywords={Hales--Jewett numbers, bracket Hales--Jewett numbers, van der Waerden numbers, power-residue colourings, CryptoMiniSat, local search}}

\newtheorem{theorem}{Theorem}
\newtheorem{proposition}[theorem]{Proposition}
\newtheorem{lemma}[theorem]{Lemma}
\theoremstyle{remark}
\newtheorem{remark}[theorem]{Remark}

\newcommand{\HJ}{\mathrm{HJ}}
\newcommand{\W}{\mathrm{W}}
\newcommand{\F}{\mathbb{F}}
\newcommand{\Z}{\mathbb{Z}}
\newcommand{\T}{\mathcal{T}}
\newcommand{\ow}{\mathrm{ow}}

\DeclareMathOperator{\ind}{ind}
\DeclareMathOperator{\type}{type}

\setlist{nosep,leftmargin=*}
\begin{document}
\title{More lower bounds for Hales--Jewett and van der Waerden numbers}
\author{Younes Mouhib\\\textnormal{Independent researcher}}
\date{17 November 2026}

\begin{abstract}
We give explicit colourings improving the comparison lower bounds assembled here for several Hales--Jewett and van der Waerden numbers. The main algebraic tool generalises Rabung's power-residue colouring from arithmetic progressions to arbitrary sets. If every translate of a $t$-element set $S\subset\F_p$ meets at least two classes of the $r$-th power character, and the integer weights have both parities, then a period-$2p$ one-weight colouring of $[t]^p$ has no monochromatic combinatorial line. The 19 displayed algebraic certificates comprise 18 improvements and one closed-form reproduction; they include $\HJ(4,3)\ge164$, $\HJ(5,4)\ge7418$, and $\HJ(5,6)\ge34034$. The underlying period-$p$ construction also proves $\kappa_{\ow}(t,r)\ge p-1$ for the one-weight bracket ceiling, in particular $\kappa_{\ow}(5,4)\ge7416$ and $\kappa_{\ow}(5,6)\ge34032$. Warm breakout-weighted tabu search on the type simplex gives $\HJ(3,4)\ge51$. Whole-colouring local search gives $\W(8,3)>522$, $\W(11,3)>1193$, and $\W(17,3)>3620$. All accepted certificates are available separately and are checked by direct enumeration; unresolved searches are reported as UNKNOWN, not UNSAT.
\end{abstract}

\maketitle

\section{Introduction}

A \emph{combinatorial line} in $[t]^n$, where $[t]=\{1,\dots,t\}$, is given by a root $\tau\in([t]\cup\{*\})^n$ with at least one $*$; its $t$ points replace every $*$ by the same letter. By the Hales--Jewett theorem~\cite{HJ}, the least $n$ for which every $r$-colouring of $[t]^n$ contains a monochromatic line exists; it is denoted $\HJ(t,r)$. Throughout, $\W(r,k)$ is the least $N$ for which every $r$-colouring of $\{1,\dots,N\}$ contains a monochromatic $k$-term arithmetic progression (AP). Thus an AP-free colouring of length $N$ proves $\W(r,k)>N$. Apart from $\HJ(2,r)=r$ and $\HJ(3,2)=4$~\cite{HT}, no nontrivial exact value is known.

With letters read as $0,\dots,t-1$, the sums of the points of a line form a $t$-term AP, so
\begin{equation}\label{eq:vdw}
\W(r,t)>N\quad\Longrightarrow\quad
\HJ(t,r)\ge\Big\lfloor\frac{N-1}{t-1}\Big\rfloor+1 .
\end{equation}
Symmetric and one-weight colourings can do substantially better~\cite{Mthesis,M1,M2}. The calibrated bounds $\HJ(3,3)\ge22$ and $\HJ(4,2)\ge14$ are reproduced in the present project, but are not claimed as new improvements.

\subsection*{Results}
\begin{itemize}
\item Table~\ref{tab:hj} gives 19 generalised-Rabung certificates: 18 improvements over the stated comparison column and the reproduction at $(4,2)$. Each $(p,g,S)$ certifies both $\HJ(t,r)\ge p+1$ and $\kappa_{\ow}(t,r)\ge p-1$. In particular, the deterministic scan at $p=34033$ gives $S=\{0,1,654,7769,23051\}$ and replaces the earlier $\HJ(5,6)\ge33224$ by $\HJ(5,6)\ge34034$.
\item A colouring of the $1326$ types of $[3]^{50}$ proves $\HJ(3,4)\ge51$ (Section~\ref{sec:hj34}).
\item Whole-colouring tabu search proves $\W(8,3)>522$, $\W(11,3)>1193$, and $\W(17,3)>3620$ (Section~\ref{sec:vdw}). By~\eqref{eq:vdw}, these give $\HJ(3,8)\ge261$, $\HJ(3,11)\ge597$, and $\HJ(3,17)\ge1810$. Baur's external certificate $\W(10,3)>938$~\cite{Baur} remains stronger than our ten-colour searches.
\end{itemize}

\begin{table}[t]
\centering\footnotesize
\renewcommand{\arraystretch}{1.10}\setlength{\tabcolsep}{2.6pt}
\caption{Generalised-Rabung certificates. Every row proves $\HJ(t,r)\ge p+1$ and $\kappa_{\ow}\ge p-1$. ``Previous'' uses the cited comparison bounds, not a fresh literature-wide record claim; (a)~the construction reproduces the bound of~\cite{M1,M2}. The certificate is $(p,g,S)$, with $\omega_a=s_a$ in~\eqref{eq:col}.}\label{tab:hj}
\begin{tabular}{@{}lrrrrrl@{}}
\toprule
$(t,r)$ & previous & $\HJ\ge$ & $\kappa_{\ow}\ge$ & $p$ & $g$ & $S$\\
\midrule
$(4,2)$ & 14\textsuperscript{a} & 14 & 12 & 13 & 2 & {\scriptsize$\{0,1,3,11\}$} \\
$(4,3)$ & 98 & 164 & 162 & 163 & 2 & {\scriptsize$\{0,1,20,21\}$} \\
$(4,4)$ & 350 & 762 & 760 & 761 & 6 & {\scriptsize$\{0,1,234,528\}$} \\
$(4,5)$ & 752 & 1\,302 & 1\,300 & 1\,301 & 2 & {\scriptsize$\{0,1,113,114\}$} \\
$(4,6)$ & 3\,260 & 3\,548 & 3\,546 & 3\,547 & 2 & {\scriptsize$\{0,1,470,3078\}$} \\
$(5,2)$ & 45 & 84 & 82 & 83 & 2 & {\scriptsize$\{0,1,2,10,75\}$} \\
$(5,3)$ & 544 & 1\,292 & 1\,290 & 1\,291 & 2 & {\scriptsize$\{0,1,34,233,1092\}$} \\
$(5,4)$ & 4\,427 & 7\,418 & 7\,416 & 7\,417 & 5 & {\scriptsize$\{0,1,105,5808,6541\}$} \\
$(5,6)$ & 24\,689 & 34\,034 & 34\,032 & 34\,033 & 7 & {\scriptsize$\{0,1,654,7769,23051\}$} \\
$(6,2)$ & 227 & 504 & 502 & 503 & 5 & {\scriptsize$\{0,1,11,12,170,345\}$} \\
$(6,3)$ & 2\,239 & 6\,200 & 6\,198 & 6\,199 & 3 & {\scriptsize$\{0,1,2262,2407,2417,5095\}$} \\
$(6,4)$ & 31\,442 & 31\,470 & 31\,468 & 31\,469 & 2 & {\scriptsize$\{0,1,1812,14285,22144,25486\}$} \\
$(7,2)$ & 618 & 1\,362 & 1\,360 & 1\,361 & 3 & {\scriptsize$\{0,1,192,231,365,412,475\}$} \\
$(7,3)$ & 8\,136 & 19\,844 & 19\,842 & 19\,843 & 19 & {\scriptsize$\{0,1,8123,14105,16463,18296,19567\}$} \\
$(8,2)$ & 1\,643 & 3\,360 & 3\,358 & 3\,359 & 11 & {\scriptsize$\{0,1,29,36,1336,1648,2772,3106\}$} \\
$(8,3)$ & 34\,058 & 62\,402 & 62\,400 & 62\,401 & 17 & {\scriptsize$\{0,1,18178,25716,36836,51121,56501,57258\}$} \\
$(9,2)$ & 5\,159 & 7\,284 & 7\,282 & 7\,283 & 2 & {\scriptsize$\{0,1,3013,3042,4155,6471,7113,7195,7254\}$} \\
$(10,2)$ & 11\,498 & 15\,662 & 15\,660 & 15\,661 & 2 & {\scriptsize$\{0,1,1315,6569,9265,9976,10851,11616,13817,14828\}$} \\
$(11,2)$ & 19\,395 & 32\,174 & 32\,172 & 32\,173 & 5 & {\scriptsize$\{0,1,6222,6665,12496,12924,13064,20611,21369,21569,24589\}$} \\
\bottomrule
\end{tabular}
\end{table}

The comparison column uses the van der Waerden bounds tabulated from~\cite{Wiki,RL,Monroe}, together with the Hales--Jewett comparisons in~\cite{Mthesis,M1,M2,CFR}. For $(5,6)$, the comparison $24\,689$ follows from $\W(6,5)>98\,754$ and~\eqref{eq:vdw}. For three colours and $t\ge9$, Kim's polychromatic-block bounds~\cite{Kim} are stronger than our algebraic rows. In particular, $\W(3,9)>4\,002\,705$ and $\W(3,10)>10\,036\,978$ give
\[
\HJ(9,3)\ge500\,339,\qquad \HJ(10,3)\ge1\,115\,220
\]
by~\eqref{eq:vdw}. We therefore omit the weaker computed algebraic bounds $193\,148$ and $592\,598$ for these two pairs from Table~\ref{tab:hj}; their valid compact certificates remain in the ancillary collection.

\section{Symmetric, one-weight, and bracket colourings}\label{sec:sym}

For $w\in[t]^n$ let $\type(w)=(\#1\text{s},\dots,\#t\text{s})$, a point of
\[
\T^{(t)}_n=\{v\in\Z_{\ge0}^t:\textstyle\sum_a v_a=n\}.
\]
If the fixed part of a line has type $u$ and the line has $k$ active coordinates, its points have types $u+ke_a$, $a\in[t]$.

\begin{lemma}[{\cite[Lemma~3.1]{Mthesis}}]\label{lem:corner}
A symmetric colouring $c=\bar c\circ\type$ of $[t]^n$ is line-free if and only if no corner tuple
$\{u+ke_a:a\in[t]\}$ with $1\le k\le n$ and $u\in\T^{(t)}_{n-k}$ is monochromatic under $\bar c$.
\end{lemma}

The simplex has $\binom{n+t-1}{t-1}$ cells and $\binom{n+t-1}{t}$ corner tuples. A valid simplex colouring lifts to an unrestricted lower-bound witness; simplex UNSAT alone would concern only symmetric colourings.

A \emph{one-weight} colouring reads only $\ell(w)=\langle\omega,\type(w)\rangle$ for some $\omega\in\Z^t$. The levels of a line are $b+k\omega_a$, where $b=\langle\omega,u\rangle$. For a colouring class $C$, let $\HJ_C^{[K]}(t,r)$ restrict the sought line to at most $K$ active coordinates, and write
\[
\kappa_C(t,r)=\max\{K\ge0:\HJ_C^{[K]}(t,r)=\infty\}.
\]
Then $\kappa_{\mathrm{sum}}\le\kappa_{\ow}\le\kappa_{\mathrm{sym}}\le\kappa_{\mathrm{all}}\le\HJ-1$. Our bracket bounds concern $\kappa_{\ow}$, not $\kappa_{\mathrm{all}}$.

\section{Generalised Rabung colourings}\label{sec:rabung}

Let $p$ be prime, $r\mid p-1$, $g$ a primitive root modulo $p$, and $\chi(x)=\ind_g(x)\bmod r$ for $x\in\F_p^\ast$. For $S=\{s_1,\dots,s_t\}\subset\F_p$ of $t\ge3$ distinct residues consider
\begin{equation*}\tag{R}
\text{for every }y\in\F_p\text{, the colours }\chi(y+s_a)\text{ with }y+s_a\ne0\text{ are not all equal.}
\end{equation*}

\begin{theorem}\label{thm:main}
Let $S$ satisfy \textup{(R)} and let $\omega\in\Z^t$ with $\omega_a\equiv s_a\pmod p$. Colour $w\in[t]^n$ by $\chi(\ell(w)\bmod p)$, arbitrarily when $\ell(w)\equiv0\pmod p$. Then no line with $1\le k\le p-1$ active coordinates is monochromatic, in any dimension. Consequently
\[
\HJ(t,r)\ge p,\qquad \kappa_{\ow}(t,r)\ge p-1 .
\]
\end{theorem}

\begin{proof}
Since $p\nmid k$, write $b\equiv ky\pmod p$. The levels are $k(y+s_a)$ modulo $p$. At most one vanishes, and the nonzero colours are $\chi(k)+\chi(y+s_a)$. By (R), these are not all equal. The displayed colouring works in every dimension, which also gives the bracket statement.
\end{proof}

\begin{proposition}\label{prop:p1}
If moreover the $\omega_a$ are not all of the same parity, then the colouring of $[t]^p$
\begin{equation}\label{eq:col}
c(w)=\psi\big(\ell(w)\bmod 2p\big),\qquad
\psi(x)=\begin{cases}
\ind_g(x\bmod p)\bmod r,& p\nmid x,\\
0,& x=0,\\
1,& x=p,
\end{cases}
\end{equation}
has no monochromatic line, so $\HJ(t,r)\ge p+1$.
\end{proposition}

\begin{proof}
Theorem~\ref{thm:main} covers $k<p$. In $[t]^p$, the only line with $k=p$ is the diagonal, whose levels $p\omega_a$ include both residues $0$ and $p$ modulo $2p$; these receive different colours.
\end{proof}

\begin{remark}\label{rem:affine}
Condition (R) is invariant under $S\mapsto\lambda S+\mu$ with $\lambda\ne0$, so one may take $S=\{0,1,s_3,\dots,s_t\}$; then $\omega=S$ has both parities. Multiplication by an $r$-th power preserves $\chi$, so full-reach scans need test only scale representatives $1,g,\dots,g^{r-1}$. Condition (R) is independent of the chosen primitive root and takes $O(pt)$ steps to check. Rabung's construction~\cite{Rabung} is the case $S=\{0,1,\dots,t-1\}$.
\end{remark}

The period-$2p$ reach stops at dimension $p$. In dimension $p+1$, take a line with $p$ active coordinates and one fixed coordinate whose level contribution $b$ is not divisible by $p$. All its levels are congruent to $b$ modulo $p$, and hence receive the same colour. Thus extending these particular palettes requires a genuinely different construction, not merely a larger simplex check.

\section{Algebraic search and certification}\label{sec:search}

\subsection*{Bitset-prefix elimination}
Store $B_c=\{x\in\F_p^\ast:\chi(x)=c\}$ as a bitset. For a prefix $S'$ of $t-1$ elements, form $I_c=\bigcap_{s\in S'}(B_c-s)$. A distinct final element $x$ violates (R) exactly when, for some colour $c$:
\begin{enumerate}[label=(\Alph*)]
\item $x\in B_c-y$ for some $y\in I_c$ (no level vanishes);
\item $x=-y$ for some $y\in I_c$ (the new level vanishes); or
\item $x\in B_c+s_0$ for some $s_0\in S'$ such that $\chi(s-s_0)=c$ for all $s\in S'\setminus\{s_0\}$ (the level of $s_0$ vanishes).
\end{enumerate}
Thus one prefix excludes every bad final element at once. The released deterministic five-element scanner enumerates normalised prefixes $\{0,1,a,b\}$ and finds the least legal $x>b$. It applies (A)--(C) with machine-word bitsets and directly rechecks every hit. At $p=34033$, $r=6$, its first hit is exactly the set in Table~\ref{tab:hj}. A generic greedy variant builds prefixes by sampling candidates and minimising $\sum_c|I_c|$; weighted local search instead replaces elements of $S$ while maintaining all translate-colour counts incrementally. A first-moment pass probability of approximately $(1-r^{1-t})^p$ per random set was used only to allocate effort, never as a proof.

\subsection*{Verification}
For this consolidation, all 22 computed algebraic records---the 19 displayed rows, two rows omitted in favour of stronger external bounds, and the superseded $p=33223$ row---were rebuilt from their primitive roots and checked at every translate: $1\,047\,828$ translates in total, with zero failures of (R). The verifier also checks primality, divisibility, primitivity of $g$, distinctness and mixed parity of $S$, and the SHA-256 hash of the complete period-$2p$ palette. The algebraic data therefore constitute complete compact certificates, not solver logs.

\section{\texorpdfstring{$\HJ(3,4)\ge51$}{HJ(3,4) >= 51}}\label{sec:hj34}

By Lemma~\ref{lem:corner}, it suffices to $4$-colour the $1326$ cells of $\T^{(3)}_{50}$ without a monochromatic corner triple; there are $22\,100$ such triples. The separate file \texttt{certificates/simplex/HJ34\_T50.json} contains the $51$ rows of digits. Row $a$ and position $b$ give the colour of $(a,b,50-a-b)$. Direct enumeration finds zero monochromatic corners, proving the claimed lower bound.

\subsection*{Actual weighted-tabu implementation}
The search maintains incremental corner-colour counts and a constant-time list of bad corners. It samples eight bad corners and selects one of largest weight, then evaluates every recolouring of its three cells by weighted conflict change. A cell's previous colour is tabu for a uniform random tenure from $12$ through $40$, with aspiration when a move beats the global best. Ties are random. If the best gain is nonnegative, all current bad-corner weights are incremented with probability $1/4$; if the best gain is positive, a random move is made with probability $1/100$. After $2\,000\,000$ steps without global improvement, the best state is perturbed and the search restarts. There is no periodic weight-smoothing step in the retained implementation.

A valid colouring of $\T_{n-1}$ shifted by $e_a$ is a warm start that already satisfies every corner inside $\{v_a\ge1\}$; the copy direction matters. Starting from the supplied $1275$-cell $\T_{49}$ certificate, the successful directional run found $\T_{50}$ after $5\,748\,043$ moves and two stagnation restarts. The original seed was not logged; the complete certificate and direct verifier, rather than an invented seed, are the reproducibility record.

\subsection*{Exact CryptoMiniSat methods}
CryptoMiniSat 5.18.0~\cite{CMS,CMSRepo} was used through its native \texttt{pycryptosat} interface, with one thread, for exact calibration and SAT encoders. For cells $i$ and colours $c$, introduce $x_{i,c}$. Add one at-least-one-colour clause per cell and, for each corner tuple and colour, one clause forbidding a monochromatic tuple. At-most-one-colour clauses are unnecessary when extracting an ordinary colouring: selecting one true colour per cell cannot create a monochromatic tuple. The original $\T_{21}$, three-colour calibration had $253$ cells, $1771$ corners, $759$ variables, and $5567$ clauses, including one harmless unit clause fixing the first cell to colour $0$.

The released encoder also supports complete first-occurrence colour symmetry breaking,
\begin{equation}\label{eq:sym}
x_{i,c}\Longrightarrow\bigvee_{h<i}x_{h,c-1}\qquad(c\ge1),
\end{equation}
as an option for new experiments; this stronger scheme was not part of the original unit-clause calibration. For one-weight instances, only realised levels need variables. In the two-colour case there is one Boolean per level and two clauses per distinct line homothet. The exact $(0,2,3,5)$, $n=13$ instance has $64$ levels and $379$ homothets and reproduces $\HJ(4,2)\ge14$.

Exact Hamming-neighbourhood repair adds at-most-one-colour clauses and a sequential-counter bound on the change literals $\neg x_{i,\bar c_i}$ around a warm colouring $\bar c$. A completed UNSAT call then excludes only that Hamming ball. It is not global simplex UNSAT and says nothing about unrestricted colourings. Budget-limited calls return UNKNOWN. No $\T_{51}$ certificate or unrestricted upper-bound conclusion is claimed.

\section{Van der Waerden numbers}\label{sec:vdw}

\subsection*{Method and results}
Published three-term certificates often come from cyclic colourings, unrolled and extended at their endpoints~\cite{Heule,HeuleRepo,Komkov}. We start from Heule's cyclic expansions of lengths $510$, $1182$, and $3546$ for $8$, $11$, and $17$ colours. Endpoint extension first gives lengths $511$, $1183$, and $3547$. At every subsequent extension, one endpoint receives a random colour and the \emph{whole} colouring is repaired; no old position is fixed.

A tabu step selects a monochromatic AP and takes a best recolouring of one of its three positions, with random tie-breaking. A changed position remains tabu for a random tenure from $2$ through $8$, with aspiration, and $2\%$ of steps are random. Independent seeds and perturbation restarts are used. The released implementation maintains AP-colour counts and the bad-AP list incrementally. Warm colourings were also tested as phase hints in unrestricted one-hot CDCL SAT, with no positions fixed; the strongest retained global runs did not use first-occurrence colour breaking and did not surpass the local-search bounds below.

\begin{table}[ht]
\centering\small
\renewcommand{\arraystretch}{1.12}
\caption{Van der Waerden bounds obtained here. The comparison column is not a fresh literature-wide record audit. The last column counts all three-term APs checked directly in the final certificate.}\label{tab:vdw}
\begin{tabular}{@{}lrlrrr@{}}
\toprule
number & previous & source of comparison & cyclic start & new & APs checked\\
\midrule
$\W(8,3)$ & $>515$ & Komkov~\cite{Komkov} & $510$ & $\mathbf{>522}$ & $67\,860$\\
$\W(11,3)$ & $>1187$ & Komkov~\cite{Komkov} & $1182$ & $\mathbf{>1193}$ & $355\,216$\\
$\W(17,3)$ & $>3561$ & recurrence of~\cite{BCT} & $3546$ & $\mathbf{>3620}$ & $3\,274\,290$\\
\bottomrule
\end{tabular}
\end{table}

The next lengths $523$, $1194$, and $3621$ were not solved within their caps and remain UNKNOWN. For 17 colours, the Blankenship--Cummings--Taranchuk recurrence~\cite{BCT}
\[
\W(r,k)>q\bigl(\W(r-\lceil r/q\rceil,k)-1\bigr),\qquad q\le k\text{ prime},
\]
with $q=3$ and the older $\W(11,3)>1187$ gives the comparison $\W(17,3)>3561$. Using our $\W(11,3)>1193$ in the recurrence gives only $3579$, still below the direct $3620$ certificate. For ten colours, we retain Baur's stronger $\W(10,3)>938$~\cite{Baur}.

The three complete strings and their full SHA-256 hashes are in \texttt{certificates/vdw/}. Direct AP enumeration found zero monochromatic progressions in each. In a random colouring, the first-moment threshold for $k=3$ is approximately $N^\ast=4r^2\ln r$, namely $532$, $1161$, and $3275$ in these cases. This explains the usefulness of global repair near the threshold but plays no role in certification. Searches for $\W(4,4)$ produced no length-$1049$ certificate; unrestricted runs were UNKNOWN, so the retained bound remains $\W(4,4)>1048$.

\section{Concluding remarks}

The effective pattern is to compress first, then combine exact and stochastic search. Algebraic palettes reduce large grids to the finite condition (R); bitset prefix elimination makes that condition efficiently searchable; warm weighted local search is effective on satisfiable simplex and interval instances; and exact SAT provides calibrated models and bounded-neighbourhood conclusions. Immediate directions are parallel prime scans for larger $t$, nonperiodic or block-symmetric colourings beyond the period-$2p$ reach, genuinely diverse $\T_{51}$ basins before further exact repair, and new structured seeds for $k\ge4$ van der Waerden instances. All bracket claims here are one-weight statements and do not settle the unrestricted Collapse Conjecture.

\subsection*{Data and code availability}
The complete certificate collection, search programs, exact CryptoMiniSat/PySAT encoders, detailed reproduction notes, and a Python standard-library verifier are available at
\begin{center}
\url{https://github.com/ysmouhib/hj-vdw-reproducibility}.
\end{center}
The arXiv ancillary directory \texttt{anc/} contains the same certificate collection and verifier. From the unpacked source, \texttt{python3 anc/verify\_all.py} checks every ancillary certificate. Nothing in the mathematical claims depends on accepting a stochastic or solver-produced witness without this direct check.

\subsection*{Acknowledgements}
The searches and drafts used AI assistance (Claude, Anthropic, and Kimi, Moonshot AI). All retained certificates were checked independently of the generating search. We thank M.~Heule and the authors of~\cite{HeuleRepo} for making their certificates public, and the CryptoMiniSat authors for their open-source solver.

\end{document}